\documentclass[12pt]{article}
\usepackage[]{amsmath,amssymb}
\usepackage{latexsym}
\usepackage{cite}

\newtheorem{definition}{Definition}[section]
\newtheorem{theorem}[definition]{Theorem}
\newtheorem{lemma}[definition]{Lemma}
\newtheorem{corollary}[definition]{Corollary}

\newtheorem{note}[definition]{Note}

\newtheorem{proposition}[definition]{Proposition}

\def\F{\mathbb F}
\def\K{\mathbb F}

\begin{document}
\title{\bf Mock Tridiagonal Systems}
\author{
Tatsuro Ito\footnote{Supported in part by JSPS grant
18340022.} $\;$   and
Paul Terwilliger\footnote{This author gratefully acknowledges 
support from the FY2007 JSPS Invitation Fellowship Program
for Reseach in Japan (Long-Term), grant L-07512.}
}
\date{}

\maketitle
\begin{abstract}
We introduce the notion of a {\it mock tridiagonal system}.
This is a generalization of a tridiagonal system
in which the irreducibility assumption
is replaced by a certain nonvanishing condition.
We show how mock tridiagonal systems
can be used to construct
tridiagonal systems that meet certain specifications.
This paper is 
part of our ongoing project to classify the
tridiagonal systems up to isomorphism.

 \end{abstract}

\section{Tridiagonal systems and mock tridiagonal systems}

\noindent 
The concept of a tridiagonal system
was introduced in
\cite[Definition~2.1]{TD00}
as a natural generalization of a Leonard system
\cite{LS99,
qrac,
madrid}
and as a tool for studying ($P$ and $Q$)-polynomial
association schemes
\cite{BanIto,
Cau,
pasc,
TersubI}.
One can view the concept as part of the
bispectral problem
\cite{GH7,GH1,Zhidd}.
There are connections to representation theory
\cite{
hasan2,
neubauer,
Ha,
tdanduq,
NN,
IT:aug,
qtet,
Ev,
IT:Krawt,
qracah,
qSerre
}
and statistical mechanical models
 \cite{bas1,bas2,bas6,bas7,DateRoan2, Davfirst, Da,
Klish1
}.
More results on tridiagonal systems can be found in
\cite{
TD00,
shape,
nom4,
N:aw,
N:refine,
nomsplit,
nomsharp,
nomtowards,
nomstructure,
nom:mu,
drin
}.
It remains an open problem to classify the
tridiagonal systems up to isomorphism, but
 classifications 
do exist for some special cases
\cite{
 IT:Krawt,
IT:aug,
NN,
Vidar,
nom:mu,
qracah
}.
To make further progress on the classification problem,
in this paper we introduce the notion of a
{\it mock tridiagonal system}.
This is a generalization of a tridiagonal system
in which the irreducibility assumption
is replaced by a certain nonvanishing condition.
In our main result, we show how mock tridiagonal systems
can be used to construct
tridiagonal systems that meet certain specifications.

\medskip
\noindent Before going into more detail we recall
the definition of a tridiagonal system.
 We will use the following terms.
Throughout this paper $\K$ denotes a field, and
 $V$ denotes a vector space over $\K$ with finite
positive dimension.
Let ${\rm End}(V)$ denote the $\F$-algebra of
all $\F$-linear transformations from $V$ to $V$.
 Given $A \in 
{\rm End}(V)$ 
and a
subspace $W \subseteq V$,
we call $W$ an
 {\it eigenspace} of $A$ whenever 
 $W\not=0$ and there exists $\theta \in \K$ such that 
$W=\lbrace v \in V \;\vert \;Av = \theta v\rbrace$;
in this case $\theta$ is the {\it eigenvalue} of
$A$ associated with $W$.
We say that $A$ is {\it diagonalizable} whenever
$V$ is spanned by the eigenspaces of $A$.
Assume $A$ is diagonalizable.
Let $\{V_i\}_{i=0}^d$ denote an ordering of the eigenspaces of $A$
and let $\{\theta_i\}_{i=0}^d$ denote the corresponding ordering of
the eigenvalues of $A$.
For $0 \leq i \leq d$ define $E_i \in 
\mbox{\rm End}(V)$ 
such that $(E_i-I)V_i=0$ and $E_iV_j=0$ for $j \neq i$ $(0 \leq j \leq d)$.
Here $I$ denotes the identity of $\mbox{\rm End}(V)$.
We call $E_i$ the {\em primitive idempotent} of
 $A$ corresponding to $V_i$ (or $\theta_i$).
Observe that
(i) $I=\sum_{i=0}^d E_i$;
(ii) $E_iE_j=\delta_{i,j}E_i$ $(0 \leq i,j \leq d)$;
(iii) $V_i=E_iV$ $(0 \leq i \leq d)$;
(iv) $A=\sum_{i=0}^d \theta_i E_i$.
Moreover
\begin{equation}         \label{eq:defEi}
  E_i=\prod_{\stackrel{0 \leq j \leq d}{j \neq i}}
          \frac{A-\theta_jI}{\theta_i-\theta_j}.
\end{equation}
Note that each of 
$\{A^i\}_{i=0}^d$,
$\{E_i\}_{i=0}^d$ is a basis for the $\K$-subalgebra
of $\mbox{\rm End}(V)$ generated by $A$.
Moreover $\prod_{i=0}^d(A-\theta_iI)=0$.

\begin{definition}
\label{def:tds}
\rm 
\cite[Definition 2.1]{TD00}
By a {\it tridiagonal system} (or {\it  TD system})
 on $V$ we mean a sequence
$\Phi=(A;\{E_i\}_{i=0}^d;A^*;\{E^*_i\}_{i=0}^\delta)$
that satisfies (i)--(iv) below.
\begin{itemize}
\item[(i)]
Each of $A,A^*$ is a diagonalizable element of $\mbox{\rm End}(V)$.
\item[(ii)]
$\{E_i\}_{i=0}^d$ is an ordering of the 
 primitive idempotents of $A$ such that
\begin{eqnarray*}
E_iA^*E_j =0 \qquad \mbox{if} \quad |i-j|>1, \qquad
(0 \leq i,j\leq d).
\end{eqnarray*} 
\item[(iii)]
$\{E^*_i\}_{i=0}^\delta$ is an ordering
of the primitive idempotents of $A^*$ such that
\begin{eqnarray*}
E^*_iAE^*_j =0 \qquad \mbox{if} \quad |i-j|>1, \qquad
(0 \leq i,j\leq \delta).
\end{eqnarray*} 
\item[(iv)] There does not exist a subspace $W$ of 
$V$ such that $AW\subseteq W$,
$A^*W \subseteq W$,
$W\not=0$, $W\not=V$.
\end{itemize}
We say that $\Phi$ is {\it over $\F$}.
\end{definition}

\begin{note}
\label{lem:convention}
\rm
According to a common notational convention $X^*$ denotes 
the conjugate-transpose of $X$. We are not using this convention.
For the TD system in Definition
\ref{def:tds} the linear transformations
$A$, $E_i$, $A^*$, $E^*_i$
are arbitrary 
subject to
(i)--(iv) above.
\end{note}

\begin{definition}
\label{def:diam}
\rm
Referring to the TD system $\Phi$ in Definition \ref{def:tds},
it turns out that
the integers $d$ and $\delta$
are equal
\cite[Lemma 4.5]{TD00}; 
 we call this common value the {\it diameter}
of $\Phi$.
\end{definition}

\noindent 
In the theory of TD systems the following situation
often occurs:
we wish to show that there exists a TD system
that meets some given specifications
 \cite{nom:mu, qracah}.
Suppose we have a candidate 
$(A;\{E_i\}_{i=0}^d;A^*;\{E^*_i\}_{i=0}^d)$,
and we just want to check that it really is
a TD system.
It is usually 
routine to 
verify 
conditions
(i)--(iii)
of
Definition 
\ref{def:tds},
but 
often difficult to verify
condition (iv) of
Definition 
\ref{def:tds}. In this paper
we give a method for
constructing TD systems
that overcomes this difficulty.
The method is based on the notion of
a mock tridiagonal system which we now introduce.

\begin{definition}
\label{def:mtds}
\rm
By a {\it mock tridiagonal system} (or {\it  MTD system})
 on $V$ we mean a sequence
$\Phi=(A;\{E_i\}_{i=0}^d;A^*;\{E^*_i\}_{i=0}^d)$
that satisfies (i)--(iv) below.
\begin{itemize}
\item[(i)]
Each of $A,A^*$ is a diagonalizable element of $\mbox{\rm End}(V)$.
\item[(ii)]
$\{E_i\}_{i=0}^d$ is an ordering of the 
 primitive idempotents of $A$ such that
\begin{eqnarray*}
E_iA^*E_j =0 \qquad \mbox{if} \quad |i-j|>1, \qquad
(0 \leq i,j\leq d).
\end{eqnarray*} 
\item[(iii)]
$\{E^*_i\}_{i=0}^d$ is an ordering
of the primitive idempotents of $A^*$ such that
\begin{eqnarray*}
E^*_iAE^*_j =0 \qquad \mbox{if} \quad |i-j|>1, \qquad
(0 \leq i,j\leq d).
\end{eqnarray*} 
\item[(iv)] Each of $E^*_0E_0E^*_0$, $E^*_0E_dE^*_0$ is nonzero.
\end{itemize}
We say that $\Phi$ is {\it over $\F$}.
\end{definition}

\begin{lemma}
\label{lem:tdismtd}
Any TD system on $V$ is an MTD system on $V$.
\end{lemma}
\noindent {\it Proof:} 
Let 
 $(A;\{E_i\}_{i=0}^d;A^*$; $\{E^*_i\}_{i=0}^d)$ 
denote the TD system in question.
We show that each of $E^*_0E_0E^*_0$,
 $E^*_0E_dE^*_0$ is nonzero.
By
\cite[Lemma~5.1]{nomsharp}
the map $E^*_0V \to E_0V$, $u \mapsto E_0u$
is  bijective and the map
 $E_0V \to E^*_0V$, $v \mapsto E^*_0v$ is bijective.
By construction $E^*_0$ acts as the identity
on $E^*_0V$.
Therefore the restriction of
$E^*_0E_0E^*_0$ to
$E^*_0V$ gives a bijection
$E^*_0V \to E^*_0V$.
The space $E^*_0V$ is nonzero so  
$E^*_0E_0E^*_0$ is nonzero.
One similarly finds that
 $E^*_0E_dE^*_0$ is nonzero.
\hfill $\Box$ \\

\noindent In Section 4 
we display some MTD systems that are not  TD systems.



\section{Statement of the main theorem}

\noindent In this section we state our main results.
In order to do this concisely we first discuss
some basic concepts.

\begin{definition}        \label{def}
\rm
Let $\Phi=(A;\{E_i\}_{i=0}^d;A^*$; $\{E^*_i\}_{i=0}^d)$ 
denote an MTD system on $V$.
For $0 \leq i \leq d$ let $\theta_i$ (resp. $\theta^*_i$)
denote the eigenvalue of $A$ (resp. $A^*$)
associated with the eigenspace $E_iV$ (resp. $E^*_iV$).
We call $\{\theta_i\}_{i=0}^d$ (resp. $\{\theta^*_i\}_{i=0}^d$)
the {\em eigenvalue sequence}
(resp. {\em dual eigenvalue sequence}) of $\Phi$.
Observe that $\{\theta_i\}_{i=0}^d$ (resp. $\{\theta^*_i\}_{i=0}^d$) are
mutually distinct and contained in $\K$.
\end{definition}
\medskip

\noindent 
Referring to Definition \ref{def}, we call $\Phi$ 
{\it sharp} whenever $E^*_0V$ has dimension 1
\cite[Definition~1.5]{nomsharp}.
By \cite[Theorem~1.3]{nomstructure}, every
TD system over an algebraically closed field is sharp.

\medskip
\noindent 
The following 
notation will be useful.

\begin{definition}
\rm
Let $\lambda$ denote an indeterminate and let $\K[\lambda]$
denote the $\K$-algebra consisting of the polynomials
in $\lambda$ that have all coefficients in $\K$.
Let $\lbrace \theta_i\rbrace_{i=0}^d$  and
 $\lbrace \theta^*_i\rbrace_{i=0}^d$
denote scalars in $\K$.
For $0 \leq i \leq d$ we define the following polynomials in
$\K[\lambda]$:
\begin{eqnarray*}
 \tau_i &=& 
  (\lambda-\theta_0)(\lambda-\theta_1)\cdots(\lambda -\theta_{i-1}), \\
 \eta_i &=&
  (\lambda-\theta_d)(\lambda-\theta_{d-1})\cdots(\lambda-\theta_{d-i+1}),  \\
 \tau^*_i &=&
  (\lambda-\theta^*_0)(\lambda-\theta^*_1)\cdots(\lambda-\theta^*_{i-1}), \\
 \eta^*_i &=&
  (\lambda-\theta^*_d)(\lambda-\theta^*_{d-1})\cdots(\lambda-\theta^*_{d-i+1}).
\end{eqnarray*}
Note that each of $\tau_i$, $\eta_i$, $\tau^*_i$, $\eta^*_i$ is
monic with degree $i$.
\end{definition}

\noindent 

\begin{lemma}
\label{lem:split}
Let 
$(A; \lbrace E_i\rbrace_{i=0}^d; A^*; \lbrace E^*_i\rbrace_{i=0}^d)$
denote a sharp MTD system on $V$, with eigenvalue
sequence $\lbrace \theta_i \rbrace_{i=0}^d$
and dual eigenvalue sequence
 $\lbrace \theta^*_i \rbrace_{i=0}^d$.
Then 
 for $0 \leq i \leq d$
there exists a unique $\zeta_i \in \K$ such that 
\begin{eqnarray*}
E^*_0 \tau_i(A) E^*_0 = 
\frac{\zeta_i E^*_0}{
(\theta^*_0-\theta^*_1) 
(\theta^*_0-\theta^*_2) 
\cdots
(\theta^*_0-\theta^*_i)}. 
\end{eqnarray*}
\end{lemma}
\noindent {\it Proof:} 
The given MTD system is sharp so $E^*_0V$ has dimension 1. Therefore
$E^*_0$ has rank 1 and this implies $E^*_0 {\mathcal A} E^*_0 = \F E^*_0$,
where ${\mathcal A}={\rm End}(V)$. The result follows.
\hfill $\Box$ \\

\begin{definition}
\label{def:split}
\rm
Let $\Phi$ denote a sharp MTD system. By the
{\it split sequence} of $\Phi$ we mean the 
sequence $\lbrace \zeta_i\rbrace_{i=0}^d$ from
Lemma
\ref{lem:split}.
Observe that $\zeta_0=1$.
\end{definition}

\begin{definition}
\label{def:pa}
\rm
Let $\Phi$
denote a sharp MTD system.
By the {\it parameter array} of $\Phi$ we mean the
sequence
 $(\{\theta_i\}_{i=0}^d; \{\theta^*_i\}_{i=0}^d; \{\zeta_i\}_{i=0}^d)$
where
 $\{\theta_i\}_{i=0}^d$
(resp. 
$\{\theta^*_i\}_{i=0}^d$)
is the eigenvalue sequence 
(resp. dual eigenvalue sequence)
of $\Phi$ and
$\{\zeta_i\}_{i=0}^d$ is the split sequence of $\Phi$.
\end{definition}

\noindent 
The following proposition
indicates the importance of the
parameter array.
The proposition
 refers to the notion of isomorphism for TD systems,
which 
is defined in
\cite[Section~3]{nomsharp}.

\begin{proposition}
\label{thm:isopa}
{\rm 
\cite[Theorem~1.6]{nomstructure}}
Two sharp TD systems over $\K$ are isomorphic
if and only if they have the same parameter array.
\end{proposition}

\noindent We now state our main result.
\begin{theorem}
\label{thm:main}
Let $\Phi$ denote a sharp MTD system over $\F$. Then
there exists a sharp TD system over $\F$ that has the
same parameter array as $\Phi$.
\end{theorem}

\noindent We will use the following strategy to
prove 
Theorem \ref{thm:main}.
Let 
$\Phi=(A; \lbrace E_i\rbrace_{i=0}^d; A^*; \lbrace E^*_i\rbrace_{i=0}^d)$
denote a sharp MTD system on $V$, with parameter array
 $(\{\theta_i\}_{i=0}^d; \{\theta^*_i\}_{i=0}^d; \{\zeta_i\}_{i=0}^d)$.
Let $T$ denote the $\F$-subalgebra of ${\rm End}(V)$ generated by
$A,A^*$, and consider the $T$-module $TE^*_0V$.
We show that $TE^*_0V$ contains a unique maximal proper
$T$-submodule. Denote this submodule by $M$ and consider
the quotient $T$-module $L=TE^*_0V/M$.
By construction the $T$-module $L$ is nonzero and irreducible.
We show that
the sequence $(A; \lbrace E_i\rbrace_{i=0}^d; A^*; \lbrace E^*_i\rbrace_{i=0}^d)$
acts on $L$ as a sharp TD system with parameter array
 $(\{\theta_i\}_{i=0}^d; \{\theta^*_i\}_{i=0}^d; \{\zeta_i\}_{i=0}^d)$.

\section{The proof of the main theorem}

\noindent In this section we give a proof of
Theorem
\ref{thm:main}.
Throughout this section we fix
a sharp MTD system 
$(A; \lbrace E_i\rbrace_{i=0}^d; A^*; \lbrace E^*_i\rbrace_{i=0}^d)$
on $V$. 

\begin{definition}
\label{def:talgebra}
\rm
Let $T$ denote the $\F$-subalgebra of 
${\rm End}(V)$ 
generated by $A,A^*$. 
By definition $T$ contains
the identity $I$ of 
${\rm End}(V)$.
By (\ref{eq:defEi}) the algebra
$T$ contains each of 
$E_i, E^*_i$ for
$0 \leq i \leq d$.
\end{definition}

\noindent We now consider the $T$-module $TE^*_0V$.
We will be discussing proper 
$T$-submodules of $TE^*_0V$.
The word {\it proper} means that the
$T$-submodule in question is properly
contained in $TE^*_0V$, or in other words
not equal to $TE^*_0V$.

\begin{definition}
\rm
Let $W$ denote a proper $T$-submodule of $TE^*_0V$.
 Then $W$ is called {\it maximal}
whenever $W$ is not contained in any proper $T$-submodule of $TE^*_0V$, 
besides itself.
\end{definition}

\noindent
Our first goal is to show that 
$TE^*_0V$ has a unique maximal proper $T$-submodule.

\begin{lemma}    \label{lem:110}  
Let $W$ denote a proper $T$-submodule of $TE^*_0V$.
 Then $E^*_0W=0$.
\end{lemma}
\noindent {\it Proof:} 
Suppose $E^*_0W \not=0$. 
The space $E^*_0V$ contains
$E^*_0W$ and has dimension 1, so
$E^*_0V=E^*_0W$. The space $W$ is $T$-invariant
and $E^*_0\in T$ so $E^*_0W \subseteq W$. Therefore
$E^*_0V \subseteq W$, which yields 
$TE^*_0V \subseteq W$. This contradicts
the fact that $W$ is properly contained in $TE^*_0V$.
Therefore 
 $E^*_0W =0$. 
\hfill $\Box$ \\

\begin{lemma}     \label{lem:ww}
Let $W$ and $W'$ denote proper $T$-submodules of $TE^*_0V$.
 Then $W+W'$ is a proper
$T$-submodule of $TE^*_0V$.
\end{lemma}
\noindent {\it Proof:} 
We show $W+W'\not=TE^*_0V$. 
Define $K=\lbrace u \in TE^*_0V\,|\,E^*_0u=0\rbrace$.
Then $K\not=TE^*_0V$, since
 $0 \not=E^*_0V \subseteq TE^*_0V$
and $E^*_0$ acts as the identity on $E^*_0V$.
The space $K$ contains each of $W,W'$ by
 Lemma \ref{lem:110}, so $K$ contains
$W+W'$.
Therefore  $W+W' \not=TE^*_0V$ 
and the result follows.
\hfill $\Box$ \\


\begin{lemma}     \label{lem:120}
There exists a unique maximal proper $T$-submodule of $TE^*_0V$.
\end{lemma}
\noindent {\it Proof:} 
Concerning existence, 
let $W$ denote the subspace of $TE^*_0V$ spanned by
all the proper $T$-submodules of $TE^*_0V$.
Then $W$ is a proper $T$-submodule of $TE^*_0V$ by
Lemma \ref{lem:ww}, and since
$TE^*_0V$ has finite dimension.
The $T$-submodule $W$ is maximal by the construction.
Concerning uniqueness, suppose $W$ and $W'$ are 
maximal proper $T$-submodules 
of $TE^*_0V$. By Lemma \ref{lem:ww} $W+W'$ is a proper
 $T$-submodule of $TE^*_0V$.
The space $W+W'$ contains each of $W$, $W'$ so $W+W'$ is equal to each of
$W$, $W'$ by the maximality of $W$ and $W'$. Therefore $W=W'$ and the result 
follows.
\hfill $\Box$ \\

\begin{definition}  
\label{def:m} 
\rm
Let $M$ denote the maximal proper $T$-submodule of $TE^*_0V$.
Let $L$ denote the quotient $T$-module $TE^*_0V\slash M$.
By construction the $T$-module
 $L$ is nonzero  and irreducible.
\end{definition}


\begin{proposition}           \label{thm:exist}  
The sequence $(A; \{E_i\}_{i=0}^d;A^*; \{E^*_i\}_{i=0}^d)$ acts on $L$ as 
a  TD system with parameter array
 $(\{\theta_i\}_{i=0}^d; \{\theta^*_i\}_{i=0}^d; \{\zeta_i\}_{i=0}^d)$.
\end{proposition}

\noindent {\it Proof:} 
The following equations hold in $T$ and hence on $L$:
\begin{eqnarray}
&&E_iE_j = \delta_{i,j}E_i,
\qquad 
E^*_iE^*_j = \delta_{i,j}E^*_i \qquad (0 \leq i,j\leq d),
\label{eq:eprod}
\\
&&  \qquad I = \sum_{i=0}^d E_i,
\quad \qquad 
I = \sum_{i=0}^d E^*_i,
\label{eq:esum}
\\
&&  \qquad A=\sum_{i=0}^d \theta_i E_i,
\qquad
\quad
A^*=\sum_{i=0}^d \theta^*_i E^*_i,
\label{eq:asum}
\\
&&E_iA^*E_j=0 \quad {\rm if} \quad |i-j|>1 \qquad \quad (0 \leq i,j\leq d),
\label{eq:ij}
\\
&&E^*_iAE^*_j=0 \quad {\rm if} \quad |i-j|>1 \qquad \quad (0 \leq i,j\leq d).
\label{eq:ijs}
\end{eqnarray}
Define
 $S =\lbrace i |0 \leq i \leq d, \;E_iL\not=0\rbrace$
and
 $S^* =\lbrace i |0 \leq i \leq d, \;E^*_iL\not=0\rbrace$.
By the equations on the left in 
(\ref{eq:eprod}),
(\ref{eq:esum})
we have $L=\sum_{i\in S}E_iL$ (direct sum).
Using the equations on the left in (\ref{eq:eprod}),
(\ref{eq:asum}) we find that 
for $i \in S$ the space $E_iL$ is an eigenspace for $A$
with eigenvalue $\theta_i$.
By these comments  
$A$ is diagonalizable on $L$ with
eigenvalues $\lbrace \theta_i\rbrace_{i\in S}$.
By
the equation on the left in
(\ref{eq:eprod}),
for $i \in S$ the element
$E_i$ acts as the identity on $E_iL$
and vanishes on $E_jL$ for $j \not=i$ $(j \in S)$.
In other words $E_i$ acts 
 on $L$ as the 
primitive idempotent of $A$ associated with $\theta_i$.
Similarly $A^*$ is diagonalizable on $L$ with
eigenvalues $\lbrace \theta^*_i\rbrace_{i\in S^*}$,
and for $i \in S^*$ the element
$E^*_i$ acts on $L$ as the 
primitive idempotent of $A^*$ associated with $\theta^*_i$.
We now show that there exist nonnegative integers $r,k$
$(r + k \leq d)$ such that
 $S =\lbrace r,r+1,\ldots, r+k\rbrace$.
The set $S$ is nonempty since $L$ is
nonzero and equal to
$\sum_{i\in S}E_iL$.
Define
$r = {\rm min}\lbrace i \,|\,i \in S\rbrace$ and
$p = {\rm max}\lbrace i \,|\,i \in S\rbrace$.
For $r+1 \leq i \leq p-1$
we have $i \in S$; otherwise
 $E_rL+\cdots + E_{i-1}L$
is a nonzero proper $T$-submodule of
$L$, contradicting the irreducibility of
the $T$-module $L$. Now
 $S =\lbrace r,r+1,\ldots, r+k\rbrace$ where
$k=p-r$.
Similarly
there exist nonnegative integers $t,k^*$
$(t + k^* \leq d)$ such that
 $S^* =\lbrace t,t+1,\ldots, t+k^*\rbrace$.
By the argument so far, the sequence
$(A; \{E_i\}_{i=r}^{r+k}; A^*; \{E^*_i\}_{i=t}^{t+k^*})$ acts on
$L$ as a TD system.
For this system
we invoke
the first sentence of
Definition \ref{def:diam} 
to get $k = k^*$.
We now show that $r=0$.
Suppose $r\not=0$.
Then $0 \not \in S$ so
 $E_0L=0$.
This implies $E_0TE^*_0V \subseteq M$
so
 $E_0E^*_0V \subseteq M$.
In this containment we apply $E^*_0$ to both sides and use $E^*_0M=0$ to
get $E^*_0E_0E^*_0 V=0$.
This contradicts
Definition \ref{def:mtds}(iv)
 so $r=0$.
Next we show that $t=0$. 
Suppose $t\not=0$. Then $0 \not\in S^*$ so $E^*_0L=0$.
This implies $E^*_0TE^*_0V \subseteq M$ so
$E^*_0V \subseteq M$.
But then $TE^*_0V \subseteq M$
since $M$ is $T$-invariant.
This contradicts the fact that $M$ is
properly contained in $TE^*_0V$, so
 $t=0$.
We now show that $k=d$.
Suppose $k \not=d$.
 Then
$d \not\in S$ so
 $E_d L=0$.
This implies $E_dTE^*_0V \subseteq M$
so
 $E_dE^*_0V \subseteq M$.
In this containment we apply $E^*_0$ to both sides and use $E^*_0M=0$ to
get $E^*_0E_dE^*_0 V=0$.
This contradicts
Definition
\ref{def:mtds}(iv)
 so $k = d$.
We have shown $(r,t,k)=(0,0,d)$,
so now $(A; \{E_i\}_{i=0}^d; A^*; \{E^*_i\}_{i=0}^d)$ acts on
$L$ as a TD system which we denote by $\Phi$.
By construction $\Phi$ has eigenvalue sequence $\{\theta_i\}_{i=0}^d$
and dual eigenvalue sequence $\{\theta^*_i\}_{i=0}^d$.
By Lemma
 \ref{lem:split}
and Definition
\ref{def:split} the equations
\begin{eqnarray*}
 E^*_0\tau_i(A)E^*_0 = 
  \frac{\zeta_i E^*_0}
       {(\theta^*_0-\theta^*_1)(\theta^*_0-\theta^*_2)
\cdots(\theta^*_0-\theta^*_i)}
     \qquad \qquad (0 \leq i \leq d)
\label{eq:spliteq}
\end{eqnarray*}
hold on $V$.  Therefore these equations
hold in $T$ and hence on $L$.
 By this and Definition \ref{def:split} the sequence 
$\{\zeta_i\}_{i=0}^d$ is the split sequence for $\Phi$.
By these comments $\Phi$ has parameter array
 $(\{\theta_i\}_{i=0}^d; \{\theta^*_i\}_{i=0}^d; \{\zeta_i\}_{i=0}^d)$
and the result follows.
\hfill $\Box$ \\

\noindent Theorem
\ref{thm:main}
is immediate from
Proposition
\ref{thm:exist}.

\medskip
\noindent We finish this section with two
corollaries of
 Theorem 
\ref{thm:main}
and Proposition
\ref{thm:exist}. 
The first corollary is about  
 the dimensions of the $\F$-vector spaces
$L$ and $M$ from Definition
\ref{def:m}.

\begin{corollary}
\label{cor:mdim}
 The following
 {\rm (i)}, {\rm (ii)} hold.
\begin{itemize}
\item[\rm (i)]
${\rm dim}(L)\geq d+1$;
\item[\rm (ii)]
${\rm dim}(M)\leq {\rm dim}(TE^*_0V) -d-1$.
\end{itemize}
\end{corollary}
\noindent {\it Proof:} 
(i) By Proposition
\ref{thm:exist}  
 there exists a TD system  on $L$ that
has diameter $d$.
\\
\noindent (ii)
Recall that $L$ is the quotient
$TE^*_0V/M$ so
${\rm dim}(L)+{\rm dim}(M)={\rm dim}(TE^*_0V)$.
The result follows from this and (i) above.
\hfill $\Box$ \\

\medskip
\noindent In the next corollary
we list some constraints satisfied by
the parameter array of a sharp MTD system.

\begin{corollary}
\label{cor:mtdconst}
Let 
 $(\{\theta_i\}_{i=0}^d; \{\theta^*_i\}_{i=0}^d; \{\zeta_i\}_{i=0}^d)$
denote the parameter array of a sharp MTD system.
Then 
{\rm (i)}--{\rm (iii)} hold below:
\begin{itemize}
\item[\rm (i)]
$\theta_i \neq \theta_j$, 
$\theta^*_i \neq \theta^*_j$ if $i \neq j$ $(0 \leq i,j \leq d)$.
\item[\rm (ii)]
$\zeta_0=1$, $\zeta_d \neq 0$, and
\begin{equation*}   
        \label{eq:ineqc}
0 \neq \sum_{i=0}^d \eta_{d-i}(\theta_0)\eta^*_{d-i}(\theta^*_0) \zeta_i.
\end{equation*}
\item[\rm (iii)]
The expressions
\begin{equation*} 
\label{eq:betaplusonec} 
\frac{\theta_{i-2}-\theta_{i+1}}{\theta_{i-1}-\theta_i},  \qquad\qquad
  \frac{\theta^*_{i-2}-\theta^*_{i+1}}{\theta^*_{i-1}-\theta^*_i}
\end{equation*}
are equal and independent of $i$ for $2 \leq i \leq d-1$.
\end{itemize}
\end{corollary}
\noindent {\it Proof:} 
Let $\Phi$
denote the sharp MTD system in question.
Then $\Phi$ satisfies condition (i) by the last
sentence of
Definition
\ref{def}.
$\Phi$ satisfies condition (ii)
by Theorem
\ref{thm:main} and since
(ii) holds for any sharp TD system
with parameter array
 $(\{\theta_i\}_{i=0}^d; \{\theta^*_i\}_{i=0}^d; \{\zeta_i\}_{i=0}^d)$
\cite[Corollary~8.3]{nomsharp}.
$\Phi$ satisfies condition (iii)
by Theorem
\ref{thm:main} and since
(iii) holds for any TD system 
with eigenvalue sequence
 $\{\theta_i\}_{i=0}^d$
and 
dual eigenvalue sequence
 $\{\theta^*_i\}_{i=0}^d$
\cite[Theorem~11.1]{TD00}.
\hfill $\Box$ \\

\section{An example}
In this section we consider a family of
sharp MTD systems over $\F$ that have diameter 2.
In \cite{Vidar} Vidar described the 
members of this family that are TD systems.
Our focus here is on the family members
that are not TD systems. For these members we find
the space $TE^*_0V$ from below Definition
\ref{def:talgebra},
and the space $M$ from Definition
\ref{def:m}.
We also describe the induced TD system on $L=TE^*_0V/M$,
from Proposition
\ref{thm:exist}.  Throughout this section we make
 use of the work of Vidar \cite[Section~9]{Vidar}.

\medskip
\noindent From now on we fix a sequence
\begin{eqnarray}
\label{eq:padiam2}
(
\lbrace \theta_i\rbrace_{i=0}^2;
\lbrace \theta^*_i\rbrace_{i=0}^2;
\lbrace \zeta_i\rbrace_{i=0}^2)
\end{eqnarray}
of scalars in $\F$ that satisfy (i), (ii) below.
\begin{itemize}
\item[\rm (i)]
$\theta_i \neq \theta_j$, 
$\theta^*_i \neq \theta^*_j$ if $i \neq j$ $(0 \leq i,j \leq 2)$.
\item[\rm (ii)]
$\zeta_0=1$, $\zeta_2 \neq 0$, and
\begin{equation}   
        \label{eq:ineqcd2}
0 \neq \zeta_2 + \zeta_1 (\theta_0-\theta_2)(\theta^*_0-\theta^*_2)+
(\theta_0-\theta_1)(\theta_0-\theta_2)
(\theta^*_0-\theta^*_1)(\theta^*_0-\theta^*_2).
\end{equation}
\end{itemize}
Our first goal is to display an MTD system over $\F$
that has parameter array
(\ref{eq:padiam2}).

\begin{definition}
\label{def:vaas}
\rm
Let $V$ denote the vector space $\F^4$ (column vectors). Define
\begin{eqnarray*}
A=\left(
\begin{array}{ c c c c }
\theta_0 & 0 &0 &0   \\
1 & \theta_1 & 0 & 0 \\
0 & 0 & \theta_1 & 0 \\
0 & 1 &  \zeta^\times_1 &\theta_2 
\end{array}
\right),
\qquad
A^*=\left(
\begin{array}{ c c c c }
\theta^*_0 & \zeta_1 & \zeta_2 &0   \\
0 & \theta^*_1 & 0 & 0 \\
0 & 0 & \theta^*_1 & 1 \\
0 & 0  &  0 &\theta^*_2 
\end{array}
\right),
\end{eqnarray*}
where
\begin{eqnarray}
\zeta^\times_1 = \zeta_1 +  
(\theta_0-\theta_1)(\theta^*_0-\theta^*_1)-
(\theta_1-\theta_2)(\theta^*_1-\theta^*_2).
\label{eq:z1times}
\end{eqnarray}
We view $A,A^* \in {\rm End}(V)$.
\end{definition}

\begin{lemma}
\label{lem:aasdiag}
The matrix $A$ (resp. $A^*$) is diagonalizable
with eigenvalues $\lbrace \theta_i \rbrace_{i=0}^2$
(resp. 
$\lbrace \theta^*_i \rbrace_{i=0}^2$).
For $0 \leq i \leq 2$ the dimension of
the eigenspace for $A$ (resp. $A^*$)
associated with $\theta_i$ (resp. $\theta^*_i$)
is 
$\binom{2}{i}$.
\end{lemma}
\noindent {\it Proof:} 
One checks that 
 $A$ has characteristic polynomial 
$(\lambda - \theta_0) 
(\lambda - \theta_1)^2 
(\lambda - \theta_2)$ and
minimal polynomial
$(\lambda - \theta_0) 
(\lambda - \theta_1) 
(\lambda - \theta_2)$.
Our assertions for $A$ follow from this.
Our assertions for $A^*$ are similary proved.
\hfill $\Box$ \\

\begin{definition}
\label{def:primid}
\rm
For $0 \leq i \leq 2$ let $E_i$ (resp. $E^*_i$) denote
the primitive idempotent of $A$ (resp. $A^*$) associated with
$\theta_i$ (resp. $\theta^*_i$).
\end{definition}

\begin{lemma}
\label{lem:edata}
We have
\begin{eqnarray*}
&&E_0=\left(
\begin{array}{ c c c c }
1 & 0 &0 &0   \\
\frac{1}{\theta_0-\theta_1} & 0 & 0 & 0 \\
0 & 0 & 0 & 0 \\
\frac{1}{(\theta_0-\theta_1)(\theta_0-\theta_2)} & 0 &  0 & 0
\end{array}
\right),
\qquad
\quad 
\quad
E^*_0=\left(
\begin{array}{ c c c c }
1 & \frac{\zeta_1}{\theta^*_0-\theta^*_1} & \frac{\zeta_2}{\theta^*_0-\theta^*_1} &\frac{\zeta_2}{(\theta^*_0-\theta^*_1)(\theta^*_0-\theta^*_2)}\\
0 & 0 & 0 & 0 \\
0 & 0 & 0 & 0 \\
0 & 0  &  0 & 0 
\end{array}
\right),
\\
&&E_1=\left(
\begin{array}{ c c c c }
0 & 0 &0 &0   \\
\frac{1}{\theta_1-\theta_0} & 1 & 0 & 0 \\
0 & 0 & 1 & 0 \\
\frac{1}{(\theta_1-\theta_0)(\theta_1-\theta_2)} & 
\frac{1}{\theta_1-\theta_2} &  \frac{\zeta^\times_1}{\theta_1-\theta_2} & 0
\end{array}
\right),
\quad
E^*_1=\left(
\begin{array}{ c c c c }
0 & \frac{\zeta_1}{\theta^*_1-\theta^*_0} & 
\frac{\zeta_2}{\theta^*_1-\theta^*_0} &
\frac{\zeta_2}{(\theta^*_1-\theta^*_0)(\theta^*_1-\theta^*_2)}\\
0 & 1 & 0 & 0 \\
0 & 0 & 1 & \frac{1}{\theta^*_1-\theta^*_2} \\
0 & 0  &  0 & 0 
\end{array}
\right),
\\
&&E_2=\left(
\begin{array}{ c c c c }
0 & 0 &0 &0   \\
0 & 0 & 0 & 0 \\
0 & 0 & 0 & 0 \\
\frac{1}{(\theta_2-\theta_0)(\theta_2-\theta_1)} & \frac{1}{\theta_2-\theta_1}
 &  \frac{\zeta^\times_1}{\theta_2-\theta_1} & 1 
\end{array}
\right),
\quad
E^*_2=\left(
\begin{array}{ c c c c }
0 &  0 & 0 
&
\frac{\zeta_2}{(\theta^*_2-\theta^*_0)(\theta^*_2-\theta^*_1)}\\
0 & 0 & 0 & 0 \\
0 & 0 & 0 & \frac{1}{\theta^*_2-\theta^*_1} \\
0 & 0  &  0 & 1 
\end{array}
\right).
\end{eqnarray*}
\end{lemma}
\noindent {\it Proof:} 
The matrices $\lbrace E_i\rbrace_{i=0}^2$ 
are obtained using
(\ref{eq:defEi}).
One similarly obtains
$\lbrace E^*_i\rbrace_{i=0}^2$.
\hfill $\Box$ \\

\begin{lemma}
\label{lem:znz}
The following
{\rm (i)}--{\rm (iii)} hold.
\begin{itemize}
\item[\rm (i)]
$E_0A^*E_2 = 0$, $E_2A^*E_0=0$;
\item[\rm (ii)]
$E^*_0AE^*_2 = 0$, $E^*_2AE^*_0=0$;
\item[\rm (iii)]
$E^*_0E_0E^*_0 \not= 0$, $E^*_0E_2E^*_0 \not=0$.
\end{itemize}
\end{lemma}
\noindent {\it Proof:} 
(i), (ii) Routine calculation using
the matrices in Definition
\ref{def:vaas} and 
Lemma
\ref{lem:edata}.
\\
\noindent (iii) The $(1,1)$-entry of $E^*_0E_0E^*_0$
is 
$(\theta_0-\theta_1)^{-1}
(\theta_0-\theta_2)^{-1}
(\theta^*_0-\theta^*_1)^{-1}
(\theta^*_0-\theta^*_2)^{-1}$
times the expression on
the right in 
(\ref{eq:ineqcd2}). This expression
is nonzero so
$E^*_0E_0E^*_0\not=0$.
The $(1,1)$-entry of $E^*_0E_2E^*_0$
is 
$(\theta_2-\theta_0)^{-1}
(\theta_2-\theta_1)^{-1}
(\theta^*_0-\theta^*_1)^{-1}
(\theta^*_0-\theta^*_2)^{-1}$
times
 $\zeta_2$.
By assumption $\zeta_2 \not=0$
so 
 $E^*_0E_2E^*_0\not=0$.
\hfill $\Box$ \\

\begin{proposition}
\label{prop:gotmtd}
The sequence 
$(A;\{E_i\}_{i=0}^2;A^*;\{E^*_i\}_{i=0}^2)$
is an MTD system on $V$ with parameter array
$(
\lbrace \theta_i\rbrace_{i=0}^2;
\lbrace \theta^*_i\rbrace_{i=0}^2;
\lbrace \zeta_i\rbrace_{i=0}^2)$.
\end{proposition}
\noindent {\it Proof:} 
Let $\Phi$ denote the sequence in question.
To show that $\Phi$ is an MTD system on $V$,
we verify the conditions (i)--(iv) of Definition
\ref{def:mtds}.
Condition (i) holds
by
Lemma \ref{lem:aasdiag}.
Condition (ii) holds
by Definition
\ref{def:primid} and
Lemma
\ref{lem:znz}(i).
Condition (iii) holds
by Definition
\ref{def:primid} and
Lemma
\ref{lem:znz}(ii).
Condition (iv) holds by
Lemma
\ref{lem:znz}(iii).
We have verified 
the conditions of
 Definition
\ref{def:mtds},
so $\Phi$
is an MTD system on $V$.
By Lemma
\ref{lem:aasdiag} 
and Definition
\ref{def:primid}, 
 $\Phi$ has
eigenvalue sequence
$\lbrace \theta_i\rbrace_{i=0}^2$
and dual eigenvalue sequence
$\lbrace \theta^*_i\rbrace_{i=0}^2$.
Using Lemma
\ref{lem:split} and Definition
\ref{def:split}
 we find that
$\Phi$ has split sequence
$\lbrace \zeta_i\rbrace_{i=0}^2$.
By these comments $\Phi$ has parameter array
$(
\lbrace \theta_i\rbrace_{i=0}^2;
\lbrace \theta^*_i\rbrace_{i=0}^2;
\lbrace \zeta_i\rbrace_{i=0}^2)$
and the result follows.
\hfill $\Box$ \\

\begin{proposition}
{\rm (Vidar  \cite[Theorem~9.1]{Vidar})}
 The following
{\rm (i)}, {\rm (ii)} are equivalent:
\begin{itemize}
\item[\rm (i)]
 The MTD system 
from
Proposition
\ref{prop:gotmtd} 
is a TD system;
\item[\rm (ii)]
$\zeta_1 \zeta^\times_1 \not= \zeta_2$.
\end{itemize}
\end{proposition}

\noindent 
In Proposition
\ref{prop:gotmtd} we displayed an MTD system
$(A;\{E_i\}_{i=0}^2;A^*;\{E^*_i\}_{i=0}^2)$.
In what follows
we consider 
the corresponding algebra $T$ from
Definition
\ref{def:talgebra}, and
the $T$-modules 
$L, M$ from
Definition \ref{def:m}.

\begin{lemma}
\label{lem:lm}
Assume
$\zeta_1 \zeta^\times_1 = \zeta_2$.
Then the following 
{\rm (i)}, {\rm (ii)} hold.
\begin{itemize}
\item[\rm (i)]
 $TE^*_0V=V$.
\item[\rm (ii)]
The $\F$-vector space
$M$ is one-dimensional and spanned by $(0,-\zeta^\times_1,1,0)^t$.
\end{itemize}
\end{lemma}
\noindent {\it Proof:} 
(i) The span of the vector $(1,0,0,0)^t$ is
$E^*_0V$,
the span of $(0,0,0,1)^t$ is 
$E_2V$, 
the span of $(0,1,0,0)^t$ is 
$(A-\theta_0 I)E^*_0V$,
and the span of 
$(0,0,1,0)^t$ is
$(A^*-\theta^*_2I)E_2V$.
Therefore 
\begin{eqnarray}
V = E^*_0V + (A-\theta_0I)E^*_0V+ (A^*-\theta^*_2I)E_2V + E_2V
\qquad (\mbox{\rm direct sum}).
\label{eq:vspan}
\end{eqnarray}
By the form of $A$,
\begin{eqnarray}
E_2V = (A-\theta_1I)(A-\theta_0 I)E^*_0V.
\label{eq:e0e2}
\end{eqnarray}
Combining 
(\ref{eq:vspan}),
(\ref{eq:e0e2}) we find $TE^*_0V=V$.
\\
\noindent (ii)
Let $W$ denote the subspace of $V$ spanned by
the vector $(0,-\zeta^\times_1,1,0)^t$. Using
Definition 
\ref{def:vaas} one checks that $(A-\theta_1I)W=0$ and
$(A^*-\theta^*_1I)W=0$.
Therefore $W$ is a $T$-submodule of $V$.
Of course $W$ is properly contained in $V$ so
$W \subseteq  M$. Consequently 
the dimension of $M$ is at least one.
The dimension of $M$ is at most one by
Corollary
\ref{cor:mdim}(ii), so the dimension of $M$ is one. 
The result follows.
\hfill $\Box$ \\

\begin{definition}
\label{def:vibasis}
\rm Assume
$\zeta_1 \zeta^\times_1 =\zeta_2$. 
By Lemma
\ref{lem:lm}(ii)
the $\F$-vector space $L$ has a basis $\lbrace v_i\rbrace_{i=0}^2$
such that
\begin{eqnarray*}
v_0 &=& (1,0,0,0)^t + M,
\\
v_1 &=& (0,1,0,0)^t + M,
\\
v_2 &=& (0,0,0,1)^t + M.
\end{eqnarray*}
\end{definition}

\begin{proposition}
 Assume
$\zeta_1 \zeta^\times_1 =\zeta_2$. 
With respect to the basis $\lbrace v_i \rbrace_{i=0}^2$
 from Definition
\ref{def:vibasis} the matrices repesenting $A$, $A^*$ are
as follows.
\begin{eqnarray*}
A:\; \left(
\begin{array}{ c c c }
\theta_0 & 0 &0   \\
1 & \theta_1  & 0 \\
0 & 1 &  \theta_2 
\end{array}
\right),
\qquad
A^*:\; \left(
\begin{array}{ c c c }
\theta^*_0 & \zeta_1  &0   \\
0 & \theta^*_1  & \zeta^\times_1  \\
0 & 0  &\theta^*_2 
\end{array}
\right).
\end{eqnarray*}
\end{proposition}
\noindent {\it Proof:} 
Routine using
Definition
\ref{def:vaas} and
Definition \ref{def:vibasis}.
\hfill $\Box$ \\

\begin{proposition}
\label{eq:eesfin}
 Assume
$\zeta_1 \zeta^\times_1 =\zeta_2$. 
With respect to the basis $\lbrace v_i \rbrace_{i=0}^2$
 from Definition
\ref{def:vibasis} the matrices repesenting $\lbrace E_i\rbrace_{i=0}^2$
and $\lbrace E^*_i\rbrace_{i=0}^2$ are as follows.
\begin{eqnarray*}
&&E_0:\;\left(
\begin{array}{ c c c }
1 & 0 &0    \\
\frac{1}{\theta_0-\theta_1} & 0  & 0 \\
\frac{1}{(\theta_0-\theta_1)(\theta_0-\theta_2)}  &  0 & 0
\end{array}
\right),
\,\; \quad 
\quad
E^*_0:\;\left(
\begin{array}{ c c c }
1 & \frac{\zeta_1}{\theta^*_0-\theta^*_1} 
&\frac{\zeta_2}{(\theta^*_0-\theta^*_1)(\theta^*_0-\theta^*_2)}\\
0  & 0 & 0 \\
0   &  0 & 0 
\end{array}
\right),
\\
&&E_1:\;\left(
\begin{array}{ c c c }
0 &0 &0   \\
\frac{1}{\theta_1-\theta_0} & 1  & 0 \\
\frac{1}{(\theta_1-\theta_0)(\theta_1-\theta_2)} & 
\frac{1}{\theta_1-\theta_2} & 0
\end{array}
\right),
\quad
E^*_1:\;\left(
\begin{array}{ c c c }
0 & \frac{\zeta_1}{\theta^*_1-\theta^*_0} & 
\frac{\zeta_2}{(\theta^*_1-\theta^*_0)(\theta^*_1-\theta^*_2)}\\
0 & 1 & \frac{\zeta^\times_1}{\theta^*_1-\theta^*_2} \\
0 & 0  &  0  
\end{array}
\right),
\\
&&E_2:\;\left(
\begin{array}{ c c c }
0 & 0 &0    \\
0 & 0 & 0 \\
\frac{1}{(\theta_2-\theta_0)(\theta_2-\theta_1)} & \frac{1}{\theta_2-\theta_1}
  & 1 
\end{array}
\right),
\quad
E^*_2:\;\left(
\begin{array}{c c c }
0 &  0  
&
\frac{\zeta_2}{(\theta^*_2-\theta^*_0)(\theta^*_2-\theta^*_1)}\\
0  & 0 & \frac{\zeta^\times_1}{\theta^*_2-\theta^*_1} \\
0 & 0  &   1 
\end{array}
\right).
\end{eqnarray*}
\end{proposition}
\noindent {\it Proof:} 
Routine using
Lemma
\ref{lem:edata}
and
Definition \ref{def:vibasis}.
\hfill $\Box$ \\

\section{Acknowledgement}
The authors thank Kazumasa Nomura, Arlene Pascasio, and Melvin Vidar
for giving this paper a close reading and offering many valuable
suggestions.

\bigskip


\noindent Tatsuro Ito \hfil\break
\noindent Division of Mathematical and Physical Sciences \hfil\break
\noindent Graduate School of Natural Science and Technology\hfil\break
\noindent Kanazawa University \hfil\break
\noindent Kakuma-machi,  Kanazawa 920-1192, Japan \hfil\break
\noindent email:  {\tt tatsuro@kenroku.kanazawa-u.ac.jp}

\bigskip

\noindent Paul Terwilliger \hfil\break
\noindent Department of Mathematics \hfil\break
\noindent University of Wisconsin \hfil\break
\noindent 480 Lincoln Drive \hfil\break
\noindent Madison, WI 53706-1388 USA \hfil\break
\noindent email: {\tt terwilli@math.wisc.edu }\hfil\break

\end{document}